\documentclass[a4paper,12pt]{article} %[fleqn]{article}

\usepackage[fleqn]{amsmath}
\usepackage{amscd,amsthm,amssymb}
\usepackage[usenames]{color}
\usepackage[all]{xypic}
\usepackage[english]{babel}
\usepackage[latin5]{inputenc}
\makeatletter
\newcommand{\specificthanks}[1]{\@fnsymbol{#1}}% Inserts a specific \thanks symbol

\usepackage{setspace}
\newtheorem{teo}{Theorem}[section]
\newtheorem{coro}[teo]{Corollary}
\newtheorem{lema}[teo]{Lemma}

\theoremstyle{definition}
\newtheorem{defi}[teo]{Definition}

\theoremstyle{remark}

\numberwithin{equation}{subsection}

\def\A{\mathsf{A}}
\def\G{\mathsf{G}}

\def\qed{\hfill $\square$}

\def\proof{\noindent \textit{Proof: }}

\def\Spec{\mathrm{Spec}\,}

\begin{document}

\title{\textsc{A remark on the invariant theory \\ of real Lie groups}}

\author{Gordillo, A.
%\thanks{Universidad de Extremadura (SPAIN)\newline  Both authors have been partially supported by Junta de Extremadura and FEDER funds.} 
\and Navarro, J. %\textsuperscript{\specificthanks{1}\,,}
%\thanks{Corresponding author. Email address: \texttt{navarrogarmendia@unex.es}\newline Department of Mathematics, Universidad de Extremadura, Avda. Elvas s/n, 06006, Badajoz, Spain.}
\and Sancho, P.}

\date{\today}

\maketitle

\begin{abstract}
We present a simple remark that assures that the invariant theory of certain real Lie groups coincides with that of the underlying affine, algebraic $\mathbb{R}$-groups. In particular, this result applies to the non-compact orthogonal or symplectic Lie groups.

%\bigskip
%\noindent \emph{Key words and phrases:} Moduli spaces, ringed spaces, jets of linear connections, differential invariants. 
%\medskip
%\noindent MSC: 58D27, 53A55
\end{abstract}

%\tableofcontents

\section*{Introduction} 
% Grupo afin da un grupo de Lie ---> Hecho estandar
% La teoria de invariantes de los grupos de Lie (clasicos, reales) no está clara
% Argumento que se da aqui 

Let $\,(E, g)\,$ be an $n$-dimensional real vector space endowed with a non-singular metric of signature $(p,q)$ and let us consider the group $\,O_{p,q}\,$ of linear isometries. To develop the theory of invariants for this group, it is usually considered in the literature as an affine algebraic $\mathbb{R}$-group, and hence algebraic techniques are employed (\cite{Goodman}, \cite{Procesi}, \cite{Weyl}).

Nevertheless, this group also has the structure of a real Lie group. Recently, both the first and second main theorems for this group $\,O_{p,q}\,$ have been widely used in the realm of natural operations in pseudo-Riemannian geometry (\cite{Jaime}, \cite{Gilkey}, \cite{Katsylo}, \cite{Navarros}) or in the construction of {\it moduli} spaces of jets (\cite{GN}, \cite{GNS}). In these settings, the point of view is that of Differential Geometry; therefore, $\,O_{p,q}\,$ is understood as a real Lie group and, in principle, it is not clear why the invariants in the ``differentiable'' sense of this non-compact and non-connected Lie group  should coincide with those invariants computed in the ``algebraic'' sense, where the action of the affine $\,\mathbb{R}$-group encodes information, not only of real points, but also of complex ones. 

For these reasons, several authors have answered this question in the particular case of the orthogonal group, by giving proofs {\it ad hoc} which confirm that the invariants of the affine $\,\mathbb{R}$-group $\,O_{p,q}\,$ coincide with those of the corresponding Lie group  (\cite{Jaime}, \cite{Gilkey}). 

The purpose of this note is to present a simple argument that generalizes this fact to any affine, algebraic $\,\mathbb{R}$-group and its corresponding real Lie group.

To be more precise, if $\,\G = \Spec A\,$ is an affine $\,\mathbb{R}$-group, its set of rational points $\,\G(\mathbb{R})\,$ has a natural structure of real Lie group. For any linear representation  $\,E\,$ of $\,\G$, the vector space of invariant vectors $\,E^\G\,$ is, in general, a subspace of the space of invariant vectors $\,E^{\G (\mathbb{R})}\,$ under the action of the Lie group of rational points. The key observation is that, under mild assumptions, the set of rational points $\, \G(\mathbb{R})\,$ is Zariski dense in $\,\G$ (Lemma \ref{PuntosRacionalesDensos}); using this, it readily follows the isomorphism $\, E^{\G} = E^{\G(\mathbb{R})} $ (Corollary \ref{Coinciden}) and, more generally, an equivalence of categories $\, \bf{Rep}(\G) \ \equiv \  \bf{Rep}_{alg}(\G(\mathbb{R}))$ (Corollary \ref{Equivalencia}).

\section{Definitions}

Let $\,G\,$ be a real Lie group. %The notions of linear representation and invariant vector in the theory of Lie groups are the following:

\begin{defi}
A {\it linear representation} of $G\,$ consists on a real vector space $E\,$ of finite dimension, together with a morphism of real Lie groups 
$$
\rho \colon G \longrightarrow \mathit{Gl}(E)\,,
$$
where $\mathit{Gl}(E)\,$ is the real Lie group of all linear automorphisms of $E\,$.

A vector $v\in V\,$ is {\it invariant under the action of $G\,$} (or $G$-invariant) if the smooth map $\,G \longrightarrow E$, $\, g \longmapsto \rho(g)(v)\,$, is the constant map $v$. 
\end{defi}

%En particular, una {\bf funcion polinomica} $E \stackrel{f}{\longrightarrow} \mathbb{R}\,$ es {\bf $G$-invariante} si cumple:
%$$
%f(g\cdot v)= f(v)\,,\,\,\,\forall g\in G\,,v\in E\,.
%$$

\medskip
Let $\,\G = \Spec A\,$ be an affine group over $\mathbb{R}$; that is, $\,A\,$ is a finitely genera\-ted Hopf $\mathbb{R}$-algebra. The corresponding notions of linear representation and invariant vector for an algebraic group are the following:

%\begin{ejem}
%No todo grupo de Lie real $G$ se puede identificar con el conjunto de puntos racionales de un grupo algebraico afin: por ejemplo el revestimiento universal de $Sl(2 , \mathbb{R})$,
%(v\'{e}ase \cite{Fulton}, p\'{a}g. 150) o el grupo de Heisenberg de un espacio vectorial con una 2-forma (v\'{e}ase \cite{Segal}, p\'{a}g. 83) no son algebraicos.
%El grupo de Heisenberg es el siguiente: sea $E$ un $\mathbb{R}$-espacio vectorial de dimensi\'{o}n par y $\omega_2$ es una dos forma cerrada y no singular, entonces el grupo: $$ \mathcal{H}(E, \omega_2) := E \times U(1) $$ donde $U(1) := \{ z \in \mathbb{C} \colon \| z \| = 1 \}$ y la ley de grupo es: $$(e_1 , z_1 ) \cdot (e_2 , z_2 ) := (e_1 + e_2\, , \ exp( i\omega_2(e_1 , e_2))  \cdot z_1 \cdot z_2 ) $$ es un grupo de Lie que no es algebraico 
%\end{ejem}

\begin{defi}
A {\it linear representation} of $\,\G = \Spec A\,$ consists on a finite dimensional $\mathbb{R}$-vector space $\,E\,$, together with a morphism of algebraic groups
$$ \rho \colon \G \longrightarrow \mathsf{Gl}(E)\,,$$
where $\mathsf{Gl}(E)=\mathrm{Spec}\,\mathbb{R}[x_{ij}]_{det}\,$ is the affine $\mathbb{R}$-group of linear automorphisms of $E\,$; that is to say, with a morphism of Hopf $\mathbb{R}$-algebras $\rho^* \colon \mathbb{R}[x_{ij}]_{det} \longrightarrow A$.

A vector $\, v \in E\,$ is {\it $\G$-invariant} if the algebraic morphism $\, \G \to E$, $g \mapsto \rho(g)(v)\,$, is the constant map $\,v$. 
\end{defi}

\medskip
Both definitions are related, because any affine $\mathbb{R}$-group $\,\G\,$ is a smooth variety and hence its set of rational points $\, \G(\mathbb{R})\,$ has a canonical structure of Lie group.

Let $\,(E, \rho)\,$ be a linear representation of $\,\G$. The algebraic morphism $\,\rho\,$ maps real points to real points, so it restricts to a linear representation of the Lie group $\G (\mathbb{R})$:
$$ \G \supset \G(\mathbb{R})  \xrightarrow{\quad \rho \quad}  \mathit{Gl}(E) \subset \mathsf{Gl}(E) \ . $$ 

Also, let us consider the vector space $\,E^\G\,$ of $\G$-invariant vectors and the vector space $\,E^{\G(\mathbb{R})}\,$ of vectors invariant under the action of the Lie group $\G(\mathbb{R})$:
\begin{align*}
E^{\G(\mathbb{R})} &= \{ \, v \in E \, \colon \, \rho(g)(v) \, = v \, , \mbox{ for all } g\in \G(\mathbb{R}) \, \} \ . 
\end{align*}

Observe the inclusion $\, E^{\G} \subseteq E^{\G(\mathbb{R})}$. In general, for example when considering Lie groups with imaginary connected components such as $\mu_3 = \Spec \mathbb{R} [x] / (x^3 - 1)\,$,
these spaces may not coincide.

\section{Rational points of affine $\mathbb{R}$-groups}

Let $\,X = \Spec A\,$ be a smooth 
affine algebraic variety over $\,\mathbb{R}$; that is, $\,A\,$ is a finitely genera\-ted $\mathbb{R}$-algebra whose local rings are regular.
We assume that $\,X\,$ has no imaginary connected components, so that each connected component has, at least, one rational point.

\medskip
Let $\, X(\mathbb{R})\,$ denote the set of rational points of $\,X$. As $\,X\,$ is smooth, $\,X(\mathbb{R})\,$ has a canonical structure of smooth manifold of the same dimension as $\,X$. This fact allows to prove the following:

%has a canonical topology: the minimal topology for which the maps $f \colon X(\mathbb{R}) \to \mathbb{R}$, $f \in A$, are continuous. 
%\begin{lemma}\label{SmoothStructure}
%The set $X(\mathbb{R})$ of rational points has a canonical structure of smooth manifold of the same dimension than $X$.
%\end{lemma}

%\begin{proof}
%Any epimorphism  $\mathbb{R}[x_1 , \ldots , x_m] \to A \to 0 $ produces a submersion of algebraic varieties $X \hookrightarrow \mathbb{R}^m$. As $X$ is smooth, the ideal of algebraic functions vanishing on $X$ is locally generated %, on a neighbourhood of a rational point $x \in X(\mathbb{R})$, 
%by polynomials with linearly independent differentials at $x$. 

%This amounts to saying that, locally, $X(\mathbb{R})$ is a smooth sub-manifold of $\mathbb{R}^m$, whose dimension is the Krull dimension of $A$. 
%\hfill $\square$
%\end{proof}

\begin{lema}\label{PuntosRacionalesDensos}
 The set of rational points $X(\mathbb{R})$ is Zariski dense in $X$.
\end{lema}

\proof We can suppose that $\,X\,$ is connected; $X\,$ being also smooth, we can assume $\,A\,$ to be an integral domain.

We will proceed by induction on the dimension of $X\,$.
%Let us argue by induction on the dimension of $X$: 
%The quid of the argument is in the first step of the induction:
 
If $X$ is one-dimensional, the zero-set of any non-zero algebraic function $\,f \in A\,$ is finite; as $\,X(\mathbb{R})\,$ is a smooth manifold of dimension one, its cardinal is infinite and hence it is Zariski dense.

Now, let us finish by proving the induction step. Given a certain function $f\in A\,$ vanishing on the set of rational points, we wish to prove that $f\,$ vanishes on the whole $n$-dimensional $X\,$. 

Let $x\in X(\mathbb{R})\,$ be a rational point and let $\mathfrak{m}_x=\{f\in A:f(x)=0\}\,$. Since $A\,$ is integral, the localization map $A\hookrightarrow A_x\,$, where $A_x\,$ is the local ring on $x\,$ ---that is, the localization of $A\,$ by the multiplicative system $A\setminus \mathfrak{m}_x\,$---, is injective, so it suffices to obtain  $f=0\,$ in $A_x\,$.

% (\cite{Matsumura}, Thm. 20.3). 

Since $A_x\,$ is a regular local ring, for any function $\,h \in\mathfrak{m}_x\subset A\,$ with non-zero differential $\, \bar{h} \in \mathfrak{m}_x / \mathfrak{m}_x^2\,$, the quotient $A_x/(h)\,$ is also regular, therefore integral, so that $(h)\,$ is a prime ideal in $A_x\,$. Moreover, this function defines a hypersurface $\,\{ h= 0\}=\mathrm{Spec}\,A/(h)\,$ that is smooth in a certain connected neighbourhood $U\cap \{h=0\}\,$ of $\,x$. As $\,f\,$ vanishes on the rational points of this neighbourhood, the induction hypothesis implies $\,f \in (h)\subset A_x\,$.

Thus, if we choose functions $\,h_1, h_2 , \ldots \in \mathfrak{m}_x\,$ with non-proportional differentials, since all the ideals $\,(h_i)\,$ are prime in $\,A_x$, we can write $ \, f \,  = \, h_1 f_1 \, = \, h_1 h_2 f_2 = \ldots\, $; hence $\, f \in \mathfrak{m}_x^i\, $ for all $\,  i \in \mathbb{N}$, and it follows that $\, f = 0\,$ in $\,A_x$.  

\qed

\begin{coro}\label{Coinciden}
Let $\, \G = \Spec A\,$ be an affine algebraic $\mathbb{R}$-group without imaginary connected components.
For any linear representation $\, (E, \rho)\,$ of $\,\G$, the inclusion $\, E^{\G} \subseteq E^{\G(\mathbb{R})}$ is in fact an isomorphism of vector spaces:
$$ E^\G \, = \, E^{\G(\mathbb{R})} \ . $$
\end{coro}
%For any linear representation $\, (E, \rho)\,$ of $\,\G$, the inclusions $\, E^{\G} \subseteq E^{\G(\mathbb{R})}$ and $\, \A^\G \, \subseteq \, \A^{\G(\mathbb{R})} $ are in fact isomorphisms of vector spaces and $\mathbb{R}$-algebras, respectively:
%$$ E^\G \, = \, E^{\G(\mathbb{R})} \qquad , \qquad \A^\G \, = \, \A^{\G(\mathbb{R})} \ . $$

\proof Let $\, v \in E^{\G(\mathbb{R})}$. The algebraic morphism
$$ \G \longrightarrow E \quad , \quad g \, \longmapsto \rho(g)(v)$$
is constant along $\, \G(\mathbb{R})$. As $\, \G(\mathbb{R})\,$ is Zariski dense on $\,\G\,$ (Lemma \ref{PuntosRacionalesDensos}), the above map is constant on all $\,\G$, and hence $\,v\,$ belongs to $\,E^\G$.

\qed

%Applying this to the linear representation of polynomial functions on $\,E$, it follows that the inclusion $\, \A^{\G} \subseteq \A^{\G (\mathbb{R})}\,$ is an isomorphism of vector spaces; hence also of $\mathbb{R}$-algebras.

%%%%%%%%%%%%%%%%%%%%%%%

\medskip

In particular, if $(E, \rho)$ is a linear representation of $\G = \Spec A$, then its algebra of polynomial functions $\,\A := S^\bullet E^* \,$ is also a linear representation of $\G$. The argument above says that the inclusion $\,\A^\G \, \subseteq \, \A^{\G(\mathbb{R})}\,$ is in fact an isomorphism of vector spaces, and hence also an isomorphism of $\mathbb{R}$-algebras.

More generally, we can consider the category  $\,\bf{Rep}(\G)\,$ of linear representations of the affine $\mathbb{R}$-group $\G = \Spec A$, and the category $\,\bf{Rep}_{alg}(\G(\mathbb{R}))\,$ of {\it algebraic} linear representations of the Lie group of its rational points  $\,\G(\mathbb{R})$.

Simple arguments, following the lines of those exposed above, allow to prove:

\begin{coro}\label{Equivalencia} Let $\, \G = \Spec A\,$ be an affine algebraic $\mathbb{R}$-group without imaginary connected components. There exists an equivalence of categories:
$$ \bf{Rep}(\G) \ \equiv \  \bf{Rep}_{alg}(\G(\mathbb{R})) \ . $$
\end{coro}

% EJEMPLO SI NO SON ALGEBRAICAS

This corollay is not true if we do not consider algebraic representations: as an example, the action of  $\mathit{Gl}_1\,$ on $\mathbb{R}^2\,$ given by:
$$
\lambda\cdot (e_1,e_2):=(e_1,e_2+(\ln |\lambda|) e_1)\, 
$$
is not semisimple: the second projection $\,p_2 \colon 0 \oplus \mathbb{R} \to \mathbb{R}$, $(0, e_2) \mapsto e_2$,  is a $\mathrm{Gl}_1$-equivariant map 
that cannot be extended to a $\mathrm{Gl}_1$-equivariant map $\,\mathbb{R}^2 \rightarrow \mathbb{R}\,$.

\subsubsection*{Example: main theorems in the category of real Lie groups}

Let $(E,g)$ be an $\mathbb{R}$-vector space of finite dimension $n$, endowed with a non-singular metric
$g$ of signature $(n_+,n_{-})$. Let $O_g$ denote the (non-compact) real Lie group of its linear isometries $(E,g) \to (E,g)$.

Our goal here is to describe the algebra $\mathsf{A}^g_m$ of polynomial functions on $m$ vectors
$$
f \colon E \times \stackrel{m}{\ldots} \times E \longrightarrow \mathbb{R}
$$ 
which are invariant under the action of $O_g$, i.e., satisfying $f(g\cdot v)= f(v)$, for all  $\, g\in G$, $\, v\in E$.
In order to do so,  let us consider the following functions $y_{ij}$, for $1 \leq i \leq j \leq m$:
$$
y_{ij} \colon E \times \stackrel{m}{\ldots} \times E \longrightarrow \mathbb{R} \quad , \quad
y_{ij} (e_1, \ldots , e_m) := g (e_i , e_j) \ .
$$

If $m > n$,  these functions have relations because the following identities hold:
\begin{equation*} %\label{Relaciones}
\left| \begin{matrix}
y_{i_0j_0} & \ldots & y_{i_0 j_n} \\
\vdots &  & \vdots \\
y_{i_nj_0} & \ldots & y_{i_nj_n}
\end{matrix}
\right| \, (e_1, \ldots , e_{m})\,   = (e_{i_0} \wedge \ldots \wedge e_{i_n} ) \cdot
( e_{j_0} \wedge \ldots \wedge e_{j_n}) = 0 \cdot 0 = 0  \ ,
\end{equation*}  where $\cdot$ denotes the metric induced by $g$ on the corresponding tensors.
%for any ordered collections $1 \leq i_0 < \ldots < i_n \leq m$, $1 \leq j_0 < \ldots < j_n \leq m$,

\begin{teo}\label{MainOrtogonal}
The algebra $\mathsf{A}_m^g$ of $O_g$-invariant polynomial functions on $E \times \stackrel{m}{\ldots}
\times E $ is generated by the functions $y_{ij}$, with $1 \leq i \leq j \leq m$.

Moreover, let $Y_{ij}$ be free symmetric variables. The map $Y_{ij}\mapsto y_{ij}$ induces a canonical isomorphism
$$\begin{CD}
  \mathbb{R} [Y_{ij}] / J_{n+1} @=  \mathsf{A}_m^g  \ ,
\end{CD}$$
where $J_{n+1}$ is the ideal generated by the functions:
\begin{equation*}\label{ThisSetting}
 J_{i_0 \ldots i_n}^{j_0 \ldots j_n} = \left| \begin{matrix}
Y_{i_0j_0} & \ldots & Y_{i_0 j_n} \\
\vdots &  & \vdots \\
Y_{i_nj_0} & \ldots & Y_{i_nj_n}
\end{matrix} \right| \ , \ \mbox{ for any } \ \ \begin{matrix}
1 \leq i_0 < \ldots < i_n \leq m\\
1\leq j_0 < \ldots < j_n \leq m
\end{matrix} \ .
\end{equation*}

In particular, if $m \leq n$, these functions $y_{ij}$ are algebraically independent.
\end{teo}

\proof The parragraph following Corollary \ref{Coinciden} allows to reduce the proof of this theorem to the corresponding statement for the algebraic orthogonal group, that is well known ({\it cfr.}, for example, \cite{Weyl} Thm. 2.17). 

\qed

\end{document}